\documentclass[a4paper]{amsart}
\title[maximal surfaces]{
 Singularities of maximal surfaces
}
\date{August 31, 2007}
\dedicatory{
 Dedicated to Yusuke Sakane on the occasion of his sixtieth birthday
}
\thanks{
 Kentaro Saji was supported by JSPS Research Fellowships
 for Young Scientists.
 Masaaki Umehara and Kotaro Yamada were supported by 
 Grant-in-Aid for Scientific Research (No.\ 15340024(B)) and
 (No. 14340024(B)), respectively.}
\usepackage[dvips]{graphicx} 
\usepackage{verbatim,enumerate}
\usepackage{amssymb}
\usepackage[usenames]{color}
\usepackage{amsthm}
\theoremstyle{plain}
 \newtheorem{theorem}{Theorem}[section]
 \newtheorem*{theorem*}{Theorem}
 \newtheorem{introtheorem}{Theorem}
   
 \newtheorem*{lemma*}{Lemma}
 \newtheorem{proposition}[theorem]{Proposition}
 \newtheorem{fact}[theorem]{Fact}
 \newtheorem*{fact*}{Fact}
 \newtheorem*{introfact}{Fact}
 \newtheorem{lemma}[theorem]{Lemma}
 \newtheorem{corollary}[theorem]{Corollary}
\theoremstyle{remark}
 
 \newtheorem{remark}[theorem]{Remark}
 \newtheorem*{remark*}{Remark}
 
 \newtheorem{example}[theorem]{Example}
\numberwithin{equation}{section}
\renewcommand{\theenumi}{{\rm(\arabic{enumi})}}
\renewcommand{\labelenumi}{\theenumi}

\newcommand{\R}{\boldsymbol{R}}
\newcommand{\C}{\boldsymbol{C}}
\newcommand{\SL}{\operatorname{SL}}

\newcommand{\SU}{\operatorname{SU}}
\newcommand{\Herm}{\operatorname{Herm}}
\newcommand{\CR}{\operatorname{CR}}
\newcommand{\CCR}{\operatorname{CCR}}
\newcommand{\CE}{\operatorname{C}}
\newcommand{\SW}{\operatorname{S}}
\newcommand{\Ker}{\operatorname{Ker}}
\newcommand{\trace}{\operatorname{trace}}
\newcommand{\rank}{\operatorname{rank}}

\newcommand{\inner}[2]{\left\langle{#1},{#2}\right\rangle}
\newcommand{\A}{\mathcal{A}}
\newcommand{\F}{\mathcal{F}}
\newcommand{\G}{\mathcal{G}}

\renewcommand{\S}{\mathcal{S}}
\renewcommand{\O}{\mathcal{O}}

\renewcommand{\Re}{\operatorname{Re}}
\renewcommand{\Im}{\operatorname{Im}}
\renewcommand{\phi}{\varphi}
\renewcommand{\epsilon}{\varepsilon}
\newcommand{\cmcone}{CMC-$1$}
\newcommand{\secondff}{\mbox{$I\!I$}}

\author{S.~Fujimori}
\address[Shoichi Fujimori]{%
   Department of Mathematics, Fukuoka University of Education,
   Munakata, Fukuoka 811-4192, Japan}
\email{fujimori@fukuoka-edu.ac.jp}

\author{K. Saji}
\address[Kentaro Saji]{%
   Department of Mathematics,
   Hokkaido University,
   Sapporo 060-0810,
   Japan
}
\email{saji@math.sci.hokudai.ac.jp}

\author{M. Umehara}
\address[Masaaki Umehara]{%
   Department of Mathematics, Graduate School of Science,
   Osaka University,
   Toyonaka, Osaka 560-0043,
   Japan
}
\email{umehara@math.sci.osaka-u.ac.jp}

\author{K. Yamada}
\address[Kotaro Yamada]{%
   Faculty of Mathematics,
   Kyushu University,
   Higashi-ku, Fukuoka 812-8581, Japan%
}
\email{kotaro@math.kyushu-u.ac.jp}

\begin{document}
\begin{abstract}
 We show that the singularities of spacelike maximal surfaces in
 Lorentz-Minkowski $3$-space generically consist of cuspidal edges,
 swallowtails and cuspidal cross caps.
 The same result holds for 
 spacelike mean curvature one surfaces in de Sitter $3$-space.
 To prove these, 
 we shall give a simple criterion for a given singular point on 
 a surface to be a cuspidal cross cap.
\end{abstract}
\maketitle
\section*{Introduction}
\begingroup
\renewcommand{\theequation}{\arabic{equation}}
In \cite{UY}, a notion of {\em maxface\/} was introduced as a
class of spacelike maximal surfaces in Lorentz-Minkowski $3$-space
with singularities (see Section~\ref{sec:maxface}). 
On a neighborhood of a singular point, a maxface
is represented by a holomorphic function using a 
Weierstrass type representation given in \cite{K}.
In fact, 
for a holomorphic function $h\in\O(U)$ defined on a simply
connected domain $U\subset\C$,
there is a maxface $f_h$ with Weierstrass data 
$(g=e^h,\omega=dz)$,
where $\O(U)$ is the set of holomorphic functions on $U$ and 
$z$ is a complex coordinate of $U$.
Conversely, for a neighborhood of a singular point of a maxface $f$,
there exists an $h\in \O(U)$ such that $f=f_h$.
For precise descriptions, see Section~\ref{sec:maxface} and \cite{UY}.

On the other hand, 
in \cite{F}, a notion of {\em \cmcone{} face\/} was introduced as a
class of spacelike surfaces of constant mean curvature one in 
de Sitter 3-space (see Section~\ref{sec:cmc1-face}),
and their global properties are investigated
(see also \cite{FRUYY}).
Like the case of maxfaces, such surfaces near a singular point
are represented by holomorphic functions; that is, 
for $h\in \O(U)$, there is a \cmcone{} face $f_h$.

In this paper, we endow the set $\O(U)$ of holomorphic functions
on $U$ with the {\em compact open $C^{\infty}$-topology}.

Then 
we shall show that cuspidal edges, swallowtails and cuspidal
cross caps are generic singularities of maxfaces in Lorentz-Minkowski 
$3$-space or \cmcone{} faces in de Sitter 3-space; that is:
\begin{introtheorem}\label{thm:maxface}
 Let $U\subset\C$ be a simply connected domain and $K$ an 
 arbitrary compact set,
 and let  $S(K)$ be the subset of  $\O(U)$ consisting of 
 $h\in \O(U)$
 such that the singular points of 
 the maxface {\rm (}resp. \cmcone{} face{\rm )} $f_h$ are cuspidal
 edges,
 swallowtails or cuspidal cross caps.
 Then $S(K)$ is an open and dense subset  of $\O(U)$.
\end{introtheorem}
We should remark that conelike singularities of maximal surfaces, 
although not generic, are still important singularities, 
which are investigated by O.~Kobayashi \cite{K}, 
Fern\'andez-L\'opez-Souam \cite{FLS} and others.
\endgroup

To prove the theroem, we shall give 
a criterion for a given singular point to be 
a cuspidal cross cap (See Theorem \ref{main}.).
It should be remarked that in \cite{KRSUY}, simple criteria for 
cuspidal edges and swallowtails have been given under
the same spirit as this paper.
These criteria in this paper and in \cite{KRSUY} 
are both really useful:
In fact, as an application, we shall show that
a duality between swallowtails and cuspidal cross caps, that
is, {\it swallowtails on maxfaces correspond to
cuspidal cross caps on their conjugate maxfaces.}
Moreover, an application of the criteria of \cite{KRSUY}
the last three authors (\cite{SUY}) 
studied the behavior of the 
Gaussian curvature near cuspidal edges and swallowtails.
Relating this result, we 
shall remark in this paper on how the behavior of the Gaussian curvature 
near a cuspidal cross cap is almost the same as that of 
a cuspidal edge (see Proposition~\ref{prop:curvature}).
It should be also remarked that, analyzing the jet spaces
of constant Gaussian curvature surfaces and 
using the criteria of \cite{KRSUY},
Ishikawa-Machida \cite{IM} showed 
{\it generic singularities on surfaces
of constant Gaussian curvature in $\R^3$ consist of
cuspidal edges and swallwtails.}

Recently, \cite{ISTa} gave simple criterions of
cuspidal lips and beaks. As a consequence,
the recognition problem of five singularities on surfaces in $\R^3$,
that is,  cuspidal edges, 
swallowtails, cuspidal cross caps,
cuspidal lips and beaks are now solved completely.
(See \cite{KRSUY},  Theorem 1.4 and \cite{ISTa}.) 

\vspace{0.5\baselineskip}

The authors thank Shyuichi Izumiya and Go-o Ishikawa for fruitful
discussions and valuable comments.
\section{Cuspidal cross caps}
\label{sec:pre}

Let $U$ be a domain in $\R^2$
and $f\colon{}U \to (N^3,g)$ a $C^{\infty}$-map 
from $U$ into a Riemannian $3$-manifold $(N^3,g)$.
The map $f$ is called a {\it frontal\/} if there exists
a unit vector field $\nu$ on $N^3$ along $f$ such that $\nu$ is
perpendicular to $f_*(TU)$. 
We call this $\nu$ is the {\em unit normal vector field} 
of a frontal $f$. 
Identifying the unit tangent bundle $T_1N^3$ with the unit
cotangent bundle $T_1^*N^3$,
the map $\nu$ is identified with the map
\[
   L=g(\nu,*):U\longrightarrow T_1^*N^3.
\]
The unit cotangent bundle $T^*_1N^3$ has a canonical
contact form $\mu$ 
and $L$ is an {\it isotropic map}, that is,
the pull back of $\mu$ by $L$ vanishes.
Namely, a frontal is the projection of an isotropic map.
We call $L$ the {\em Legendrian lift\/} 
(or {\em isotropic lift}) of $f$. 
If $L$ is an immersion, the projection $f$ is called 
a {\em front}.
Whitney \cite{Whitney} proved that the generic singularities of
$C^\infty$-maps of $2$-manifolds into $3$-manifolds 
can only be cross
caps. 
(For example, $f_{\CR}(u,v)=(u^2,v,uv)$ gives a cross cap.)
On the other hand, a cross cap is not a frontal, and 
it is also well-known that cuspidal edges and swallowtails are generic
singularities of fronts 
(see, for example,  \cite{AGV}, Section 21.6, page 336). 
The typical examples of a cuspidal edge $f_{\CE}$ and 
a swallowtail $f_{\SW}$ are given by
\[
  f_{\CE}(u,v):=(u^2,u^3,v),\qquad 
  f_{\SW}(u,v):=(3u^4+u^2v,4u^3+2uv,v).
\]
\begin{figure}[t]
\begin{center}
 \includegraphics[width=3cm]{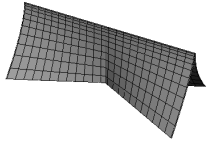}
\end{center}
\caption{The cuspidal cross cap}\label{fig:ccr}
\end{figure}
A cuspidal cross cap is a singular point which is
${\A}$-equivalent to the $C^\infty$-map 
(see Figure~\ref{fig:ccr})
\begin{equation}\label{eq:ccr}
  f_{\CCR}(u,v):= (u,v^2,uv^3),
\end{equation}
which is not a front but a frontal with unit normal vector field
\[
  \nu_{\CCR}:=\frac{1}{\sqrt{4+9u^2v^2+4v^6}}
                          (-2v^3,-3uv,2).
\]
Here, two $C^{\infty}$-maps $f\colon{}(U,p)\to N^3$ and 
$g\colon{}(V,q)\to N^3$ are {\em ${\A}$-equivalent\/} 
(or {\it right-left equivalent}) at the
points $p\in U$ and $q\in V$ if there 
exists a local diffeomorphism $\varphi$ of $\R^2$ with $\varphi(p)=q$ 
and a local diffeomorphism $\Phi$ of $N^3$ with $\Phi(f(p))=g(q)$
such that $g=\Phi\circ f \circ \varphi^{-1}$.

In this section, we shall give a simple criterion for a given singular
point on the surface to be a cuspidal cross cap.
Let  $(N^3,g)$ be a Riemannian $3$-manifold and $\Omega$ the Riemannian
volume element on $N^3$. 
Let $f:U\to N^3$ be a frontal defined on a domain $U$ in $\R^2$.
Then we can take the unit normal vector field 
$\nu\colon{}U\to T_1N^3$ of $f$ as mentioned above.
The smooth function $\lambda\colon{}U\to \R$ defined by
\begin{equation}\label{eq:lambda}
   \lambda(u,v):=\Omega(f_u,f_v,\nu)
\end{equation}
is called {\it the signed area density function},
where $(u,v)$ is a local coordinate system of $U$. 
The singular points of $f$ are the zeros of $\lambda$.
A singular point $p\in U$ is called {\it non-degenerate\/}
if the exterior derivative $d\lambda$ does not vanish at $p$.
(It can be easily checked that this non-degeneracy condition is
independent of the choice of a local coordinate system.)

When $p$ is a non-degenerate singular point, the singular set
$\{\lambda=0\}$ consists of a regular curve near $p$,
called the {\it singular curve}, 
and we can express it as a parametrized curve
$\gamma(t):(-\varepsilon,\varepsilon)\to U$
such that $\gamma(0)=p$ and
\[
  \lambda\bigl(\gamma(t)\bigr)=0 \qquad 
    (t\in (-\varepsilon,\varepsilon)).
\]
We call the tangential direction $\gamma'(t)$ 
the {\em singular direction}.
Since $d\lambda\ne 0$,  $f_u$ and $f_v$ do not vanish simultaneously. 
So the kernel of $df$ is $1$-dimensional at each singular point $p$.
A nonzero tangential vector $\eta\in T_pU$ belonging to the
kernel is called the {\it null direction}.
There exists a smooth vector field $\eta(t)$ along the singular curve
$\gamma(t)$
such that $\eta(t)$ is the null direction at $\gamma(t)$ for each $t$.
We call it the {\it vector field of the null direction}.
In \cite{KRSUY}, the following criteria for cuspidal edges and
swallowtails are given:
\begin{introfact}
 Let $f\colon{}U\to N^3$ be a front and $p\in U$ a non-degenerate
 singular point.
 Take a singular curve $\gamma(t)$ with $\gamma(0)=0$ and a vector 
 field of null directions $\eta(t)$.
 Then
 \begin{enumerate} 
  \item The germ of $f$ at $p=\gamma(0)$ is ${\A}$-equivalent to 
    a cuspidal edge if and only if
    the null direction $\eta(0)$ is transversal to the 
    singular direction $\gamma'(0)$.
  \item The germ of $f$ at $p=\gamma(0)$ is ${\A}$-equivalent to 
    a swallowtail if and only if
    the null direction $\eta(0)$ is proportional to the 
    singular direction $\gamma'(0)$ and
    \[
    \left. \frac{d}{dt}\right |_{t=0} 
                \det\bigl(\eta(t),\gamma'(t)\bigr)\ne 0,
    \]
    where
    $\eta(t)$ and $\gamma'(t)\in T_{\gamma(t)}U$ are considered
    as column vectors, and $\det$ denotes the determinant of 
    a $2\times2$-matrix.
 \end{enumerate}
\end{introfact}

We denote by $(T^*N^3)^\circ$ the complement of the zero section in
$T^*N^3$.

\begin{lemma}\label{lem:lift}
 Let $f:U\to (N^3,g)$ be a frontal. 
 Then there exists a $C^\infty$-section
 $L\colon{}U\to (T^*N^3)^\circ$  along $f$ such that 
$(\pi\circ L)_*(T_pU)\subset\Ker L_p$ for all $p\in U$,
 where $\pi\colon{}(T^*N^3)^{\circ}\to N^3$ is the canonical
 projection,
 and $\Ker L_p\subset T_{f(p)}N^3$ is the kernel of
 $L_p\colon{}T_pN^3\to \R$. 
 We shall call such a map $L$ the {\it admissible lift\/} of $f$.
 Conversely, let $L:U\to (T^*N^3)^\circ$ be a smooth section
 satisfying $(\pi\circ L)_*(T_pU)\subset \Ker L_p$.
 Then $f:=\pi\circ L$ is a frontal and $L$ is a lift of $f$.
\end{lemma}

By this lemma, we know that
the concept of frontal does not depend on the
Riemannian metric of $N^3$.
Frontal can be interpreted as a projection of
a mapping $L$ into $N^3$  satisfying 
$(\pi\circ L)_*(T_pU)\subset \Ker L_p$
for all $p\in U$.
(The projection of such an $L$ into the unit cotangent bundle
$T_1^*N^3$ gives the Legendrian lift of $f$. 
An admissible lift of $f$ is not uniquely determined,
since  multiplication of $L$ by non-constant functions also
gives 
admissible lifts.)

\begin{proof}[Proof of Lemma~\ref{lem:lift}]
 Let $\nu$ be the unit normal vector field of $f$.
 Then the map
 \[
   L:U\ni p \longmapsto g_p(\nu,*)\in T^*N^3
 \]
 gives an admissible lift of $f$. 
 Conversely, let $L\colon{}U\to (T^*N^3)^\circ$ be a non-vanishing smooth
 section with $(\pi\circ L)_*(T_pU)\subset \Ker L_p$.
 Then a non-vanishing section of 
 the orthogonal complement
 $(\Ker L)^{\perp}$ gives a normal vector field of $f$.
\end{proof}

Let $TN^3|_{f(U)}$ be the restriction of the tangent bundle $TN^3$ 
to $f(U)$. The subbundle of $TN^3|_{f(U)}$ perpendicular to
the unit normal vector $\nu$ is called the {\it limiting tangent bundle}.

As pointed out in \cite{SUY}, 
the non-degeneracy of the singular points is also independent
of the  Riemannian metric $g$ of $N^3$.
In fact, Proposition~1.3 in \cite{SUY} can be
proved under the weaker assumption that $f$ is only a frontal. 
In particular, we can show the following:
\begin{proposition}\label{non-deg}
 Let $f\colon{}U\to N^3$ be a frontal and $p\in U$ a singular point
 of $f$. 
 Let $\Omega$ be a nowhere vanishing $3$-form on $N^3$
 and $E$ a vector field on $N^{3}$ along $f$ which is transversal
 to the limiting tangent bundle.
 Then $p$ is a non-degenerate singular point of $f$ if and only if
 $\lambda=\Omega(f_u,f_v,E)$ satisfies $d\lambda\ne 0$.
\end{proposition}
Before describing the criterion for cuspidal cross caps,
we shall recall the covariant derivative along a map.
Let $D$ be an arbitrarily fixed linear connection on $N^3$
and $f:U\to N^3$ a $C^{\infty}$-map.
We take a local coordinate system $(V;x^1,x^2,x^3)$
on $N^3$ and write the connection as
\[
  D_{\frac{\partial}{\partial x^i}}\frac{\partial}{\partial x^j}
    =\sum_{k=1}^{3}\Gamma_{ij}^k \frac{\partial}{\partial x^k}.
\]
We assume that $f(U)\subset V$.
Let $X:U\to TN^3$ be 
an arbitrary vector field of $N^{3}$ along $f$
given by
\[
  X=\xi^1(u,v)\left(\frac{\partial}{\partial x^1}\right)_{f(u,v)}+
    \xi^2(u,v)\left(\frac{\partial}{\partial x^2}\right)_{f(u,v)}+
    \xi^3(u,v)\left(\frac{\partial}{\partial x^3}\right)_{f(u,v)}.
\]
Then its covariant derivative along $f$ is defined by
\[
  D^{f}_{\frac{\partial}{\partial u^l}}X:=
    \sum_{k=1}^3 \left(\frac{\partial \xi^k}{\partial u^l}
    +\sum_{i,j=1}^3 
    (\Gamma_{ij}^k\circ f) \, \xi^j \frac{\partial f^i}
    {\partial u^l}\right) \frac{\partial}{\partial x^k},
    \qquad (l=1,2),
\]
where $(u^1,u^2)=(u,v)$ is the coordinate system of $U$ 
and $f=(f^1,f^2,f^3)$.
Let 
\[
   \eta=\eta^1\frac{\partial}{\partial u^1}+\eta^1
        \frac{\partial}{\partial u^2} \in TU
\]
be a null vector of $f$, that is, $f_*\eta=0$.
In this case, we have
$\eta^1  f^k_u+ \eta^2 f^k_v=0$
for $k=1,2,3,$ and thus 
\[
  D^f_{\eta}X=
       \sum_{k=1}^3 \left(\sum_{l=1}^2 \eta^l
       \frac{\partial \xi^k}{\partial u^l}\right)
       \frac{\partial}{\partial x^k}
\]
holds, which implies the following:
\begin{lemma}\label{lem:null}
 The derivative $D_{\eta}^f X$ does not depend on the
 choice of the linear connection $D$ if $\eta$ is a
 null vector of $f$.
\end{lemma}
The purpose of this section is to prove the following assertion:
\begin{theorem}\label{main}
 Let $f\colon{}U\to N^3$ be a frontal
 and $L\colon{}U\to (T^*N^3)^\circ$  an admissible lift of $f$.
 Let $D$ be an arbitrary linear connection on $N^3$.
 Suppose that $\gamma(t)$ $(|t|<\varepsilon)$ 
 is a singular curve on $U$ passing through
 a non-degenerate singular point $p=\gamma(0)$, and that
 $X\colon{}(-\varepsilon,\varepsilon) \to TN^3$ is 
 an arbitrarily fixed vector field along $\gamma$ such that
 \begin{enumerate}
  \item\label{ass:main:1} 
       $L(X)$ vanishes on $U$, and 
  \item\label{ass:main:2} 
       $X$ is transversal to the subspace $f_*\left(T_pU\right)$
       at $p$.
 \end{enumerate}
 We set
 \[
   \tilde \psi(t):=L\bigl(D^f_{\eta(t)} X_{\gamma(t)}\bigr),
 \]
 where $\gamma(t)$ is the singular curve at $p$, $\eta(t)$
 is a null vector field along $\gamma$.
 Then the germ of $f$ at $p=\gamma(0)$ is ${\A}$-equivalent to 
 a cuspidal cross cap if and only if 
 \begin{enumerate}
 \renewcommand{\theenumi}{{\rm{(\roman{enumi})}}}
 \let\labelenumi=\theenumi
  \item\label{cond:main:1} 
       $\eta(0)$ is transversal to $\gamma'(0)$, and
  \item\label{cond:main:2} 
       $\tilde \psi(0)=0$ and $\tilde \psi'(0)\ne 0$
 \end{enumerate}
 hold, where $'=d/dt$.
\end{theorem}

\begin{remark*}
 This criterion for cuspidal cross caps  is 
 independent of the  metric and the linear connection
 of the ambient space. 
 This property will play a crucial role in Section~\ref{sec:cmc1-face}, 
 where we investigate singular points on constant mean curvature one surfaces
 in de Sitter $3$-space.
\end{remark*}

As a corollary, we shall prove the following: 

\begin{corollary}\label{thm:criterion}
 Let $f\colon{}U\to (N^3,g)$ be a frontal with unit normal vector
 field $\nu$ with the Riemannian volume element $\Omega$ on $N^3$, 
 and $\gamma(t)$ a singular curve on $U$ passing through
 a non-degenerate singular point $p=\gamma(0)$.
 We set
 \[
   \psi(t):=\Omega(\tilde \gamma',D^{f}_\eta\nu,\nu),
 \]
 where $\tilde\gamma = f\circ\gamma$,
 $D^{f}_\eta\nu$ is the canonical covariant derivative
 along a map $f$ induced from the Levi-Civita connection on $(N^3,g)$,
 and $'=d/dt$.
 Then the germ of $f$ at $p=\gamma(0)$ is ${\A}$-equivalent to 
 a cuspidal cross cap if and only if
 \begin{enumerate}
  \item[{\rm (i)}]\label{item:criterion:1}
    $\eta(0)$ is transversal to $\gamma'(0)$,
  \item[{\rm (ii)}]\label{itm:criterion:2}
    $\psi(0)=0$ and $\psi'(0)\ne 0$.
 \end{enumerate}
\end{corollary}

We remark that an application of this criterion is given in 
Izumiya-Saji-Takeuchi 
\cite{ISTe}.

\begin{proof}[Proof of Corollary~\ref{thm:criterion}]
 We set
 $X_0:=\tilde \gamma'\times_g \nu$,
 where $\tilde\gamma=f\circ\gamma$, $'=d/dt$ and 
 $\times_g$ is the vector product of $TN^3$ with respect to the
 Riemannian metric $g$.
 Since $X_0$ is perpendicular to $\nu$, we have $L(X_0)=0$.
 Moreover, $X_0$ is obviously transversal to $\tilde\gamma'$,
 and then it satisfies the conditions
 \ref{ass:main:1} and \ref{ass:main:2} in Theorem~\ref{main}. 
 On the other hand, 
 $L:=g(\nu,*)$ gives an admissible lift of $f$ and we have 
 \begin{align}\label{eq:psi}
  \psi(t)
  &=\Omega(\tilde \gamma',D^f_\eta\nu,\nu) 
  =g(\nu \times_g \tilde \gamma',D^f_\eta\nu) \\
  &=-g(X_0,D^f_\eta\nu) =g(D^f_\eta X_0,\nu) 
  =L(D^f_\eta X_0) =\tilde \psi(t). \nonumber
 \end{align}
 This proves the assertion. 
\end{proof}
To prove Theorem~\ref{main}, we prepare two lemmas.
The first one is obvious because the limiting tangent bundle
is defined as the orthogonal complement of $\tilde\nu$,
where $\tilde\nu$ is a (not necessarily unit) normal vector field.
\begin{lemma}\label{lem:front-2}
 Let $f\colon{}U \to (N^3,g)$ be a frontal and $p\in U$
 a singular point such that  $\rank (df)_p=1$ 
 {\rm (}that is, $p$ is a corank-one singular point{\rm)}.
 Take  a {\rm (}not necessarily unit{\rm )}
 normal vector field $\tilde\nu$ of $f$, that is, the limiting tangent 
 bundle is defined as $\{X\,;\,g(X,\tilde\nu)=0\}$.
 Then $f$ is a front on a neighborhood of corank-one singular point
 $p$ if and only if $D^f_{\eta}\tilde\nu$ is not proportional to 
 $\tilde\nu$ at $p$.
\end{lemma}
\begin{remark*}
 Consider a frontal $f\colon{}U\to (N^3,g)$ in
 a pseudo-Riemannian manifold with an {\em indefinite} metric $g$.
 In this case, the unit normal vector may not be defined, but there
 exists a vector field $\tilde\nu$ along $f$ such that
 the limiting tangent bundle is obtained as
 $\{X\,;\,g(\tilde\nu,X)=0\}$.
 Then  Lemma~\ref{lem:front-2} holds for $\tilde\nu$ in this case.
\end{remark*}
\begin{corollary}\label{cor:X0}
 Under the same assumptions as in Theorem~\ref{main} 
 and  the null direction is not proportional to the singular
 direction.
 Then 
 $\tilde \psi(0)\ne 0$ holds if and only if 
 $f$ is a front on a sufficiently small neighborhood of $p$.
\end{corollary}
\begin{proof}
 Take the unit normal vector field $\nu$ and the Levi-Civita connection
 $D$.
 Then $L(D^f_{\eta}X)=g(D_\eta^f X,\nu)=-g(X,D^f_{\eta}\nu)$.
 Here, $D_{\eta}\nu$ is perpendicular to $\nu$ and $df(\dot\gamma)$
 under the assumptions, we have the conclusion.
 Since $D^{f}_{\eta}\nu$ is perpendicular to both $\nu$ and
 $df(\gamma')$,
 Lemma~\ref{lem:front-2} implies the conclusion.
\end{proof}

\begin{lemma}\label{lem:X}
 Let $f\colon{}U\to N^3$ be a frontal and 
 $p$ a non-degenerate singular point of $f$ 
 satisfying  \ref{cond:main:1} of Theorem~\ref{main}.
 Then the condition $\tilde \psi(0)=\tilde \psi'(0)=0$ 
 is independent of the choice of vector field $X$ 
 along $f$ satisfying \ref{ass:main:1} and \ref{ass:main:2}.
\end{lemma}

\begin{proof}
 By \ref{cond:main:1}, we may assume that
 the null vector field $\eta(t)$ ($|t|<\varepsilon$) is 
 transversal to $\gamma'(t)$.
 Then we may take a coordinate system $(u,v)$ with the origin 
 at $p$ such that the $u$-axis corresponds to the singular curve
 and 
 $\eta(u)=(\partial/\partial v)|_{(u,0)}$.
 We fix an arbitrary vector field $X_0$ satisfying \ref{ass:main:1} 
 and \ref{ass:main:2}.
 By \ref{ass:main:2}, $X_0$ is transversal to the vector field
 $V:=f_*(\partial/\partial u)(\neq 0)$ along $f$.
 Take an arbitrary vector field $X$ along $f$ satisfying 
 \ref{ass:main:1} and \ref{ass:main:2}.
 Then it can be expressed as a linear combination
 \[
      X=a(u,v) X_0 + b(u,v)V \qquad \bigl(a(0,0)\ne 0\bigr).
 \]
 Then we have
 \[
      D^f_\eta X=da(\eta)X_0+db(\eta)V+
                     aD^f_\eta X_0 +bD^f_\eta V.
 \]
 Now 
 $L(V)=0$ holds, since $L$ is an admissible lift of $f$.
 Moreover, \ref{ass:main:1} implies that $L(X)=0$, and we have
 \[
    L(D^f_\eta X)=aL(D^f_\eta X_0) +b L(D^f_\eta V).
 \]
 Since $D^{f}_\eta$ does not depend on the choice of a connection $D$,
 we may assume that $D$ is a torsion-free connection.
 Then we have
 \[
   D^f_\eta V=
       D^f_{\frac{\partial}{\partial v}}
          f_*\left(\frac{\partial}{\partial u}\right)
   =
    D^f_{\frac{\partial}{\partial u}}
        f_*\left(\frac{\partial}{\partial v}\right)
   =0,
 \]
 since $f_*(\partial/\partial v)=f_*\eta=0$.
 Thus we have
 \[
         L(D^f_\eta X)=aL(D^f_\eta X_0)=a\tilde \psi(u).
 \]
 Since $a(0,0)\ne 0$, the conditions \ref{ass:main:1} and 
 \ref{ass:main:2}  for $X$ are the same as those of $X_0$.
\end{proof}

The following two lemmas are well-known (see \cite{GG}).
They plays a crucial role 
in Whitney \cite{Whitney} to give 
a criterion for a given $C^\infty$-map to be a cross cap. 
Let $h(u,v)$ be a $C^\infty$-function defined around the
origin.

\begin{fact}[Division Lemma]
\label{fact:division}
 If $h(u,0)$ vanishes for sufficiently small $u$,
 then there exists a $C^\infty$-function 
 $\tilde{h}(u,v)$ defined around the origin such that 
 $h(u,v)=v\tilde h(u,v)$ holds.
\end{fact}
\begin{fact}[Whitney Lemma]
\label{fact:whitney}
 If $h(u,v)=h(-u,v)$ holds for sufficiently small $(u,v)$,
 then there exists a $C^\infty$-function 
$\tilde h(u,v)$ defined around the
 origin such that $h(u,v)=\tilde h(u^2,v)$ holds.
\end{fact}
\begin{proof}[Proof of Theorem~\ref{main}]
As the assertion is local in nature, we may assume that
 $N^3=\R^3$ and let $g_0$ be the canonical metric.
 We denote the inner product associated with $g_0$ by
 $\inner{~}{~}$.
 The canonical volume form $\Omega$ is nothing but the determinant:
 $\Omega(X,Y,Z)=\det(X,Y,Z)$.
 Then the signed area density function $\lambda$ defined in
 \eqref{eq:lambda} in the introduction is written as
\[
    \lambda(u,v)= \det(f_u,f_v,\nu).
\]
 Let $f\colon{}U\to \R^3$ be a frontal and
 $\nu$  the unit normal vector field of $f$.
 Take a coordinate system $(u,v)$ centered at the singular point $p$
 such that the $u$-axis is a singular curve and the vector
 field $\partial/\partial v$ gives the null direction along the
 $u$-axis.
 If we set $X=V\times_g\nu$,
 then it satisfies \ref{ass:main:1} and \ref{ass:main:2} of
 Theorem~\ref{main}, and by \eqref{eq:psi} we have
 \[
   \tilde \psi(u)
    =\det(\gamma',D^{f}_\eta\nu,\nu)
    =\det(f_u,\nu_v,\nu).
 \]
 Thus we now suppose that $\tilde \psi(0)=0$ and 
 $\tilde \psi'(0)\ne 0$. 
 It is sufficient to show that $f$ is ${\A}$-equivalent to the standard
 cuspidal cross cap as in \eqref{eq:ccr} in the introduction.

 Without loss of generality, we may set $f(0,0)=(0,0,0)$.
 Since $f$ satisfies \ref{ass:main:1}, 
 $f(u,0)$ is a regular space curve.
 Since $f_u(u,0)\neq 0$, we may assume $f^1_u(u,0)\neq 0$ for 
 sufficiently small $u$, where we set $f=(f^1,f^2,f^3)$.
 Then the map
 \[
    \Phi\colon{}(y^1,y^2,y^3)\longmapsto 
                \bigl(f^1(y^1,0),f^2(y^1,0)+y^2,f^3(y^1,0)+y^3\bigr)
 \]
 is a local diffeomorphism of $\R^3$ at the origin.
 Replacing $f$ by  $\Phi^{-1}\circ f(u,v)$,
 we may assume
 $f(u,v) =  \bigl(u, f^2(u,v),f^3(u,v)\bigr)$,
 where  $f^2$ and $f^3$ are smooth functions around the origin
 such that $f^2(u,0)=f^3(u,0)=0$ for sufficiently small $u$.

 Then by the division lemma (Fact~\ref{fact:division}),
 there exist $C^\infty$-functions $\tilde f^2(u,v)$, $\tilde f^3(u,v)$
 such that
 $f^j(u,v)=v \tilde f^j(u,v)$ ($j=2,3$).
 Moreover, since $f_v=0$ along the $u$-axis,
 we have $\tilde f^{2}(u,0)=\tilde f^{3}(u,0)=0$.
 Applying the division lemma again, 
 there exist $C^\infty$-functions $a(u,v)$, and $b(u,v)$
 such that
 \[
    f(u,v)=\bigl(u, v^2 a(u,v),v^2 b(u,v)\bigr).
 \]
 Since $f_v(u,0)=0$, $\lambda_u(u,0)=0$ and $d\lambda\ne 0$, 
 we have
 \[
    0\ne \lambda_v(u,0)=
           \det(f_{uv},f_v,\nu)+
                \det(f_{u},f_{vv},\nu)+
                \det(f_{u},f_v,\nu_v)
                      =\det(f_{u},f_{vv},\nu).
 \]
 In particular, we have 
 \[
         0\ne f_{vv}(0,0)=2\bigl(0, a(0,0),b(0,0)\bigr).
 \]
 Hence, changing the  $y$-coordinate to the $z$-coordinate
  if necessary,
 we may assume that $a(0,0)\ne 0$.
 Then the map 
 $(u,v) \mapsto(\tilde u,\tilde v)=(u,v\sqrt{a(u,v)})$
 defined near the origin
 gives a  new local coordinate around $(0,0)$ by the
 inverse function theorem.
 Thus we may assume that
 $a(u,v)=1$, namely
 \[
          f(u,v)=\bigl(u, v^2 ,v^2 b(u,v)\bigr).
 \]
 Now we set
 \[
    \alpha(u,v):=\frac{b(u,v)+b(u,-v)}{2},\qquad
    \beta(u,v):=\frac{b(u,v)-b(u,-v)}{2}.
 \]
 Then $b=\alpha+\beta$ holds, and 
 $\alpha$ (resp.\ $\beta$) is an even (resp.\ odd) function.
 By applying the Whitney lemma, there exist
 smooth functions $\tilde\alpha(u,v)$ and $\tilde\beta(u,v)$
 such that
 \[
           \alpha(u,v)=\tilde\alpha(u,v^2),\qquad 
           \beta(u,v) =v\tilde\beta(u,v^2).
 \]
 Then we have
 \[
    f(u,v)=
     \bigl(u,v^2,v^2 \tilde\alpha(u,v^2)+v^3 \tilde\beta(u,v^2)\bigr).
 \]
 Here,
 \[
       \Phi_1:(x,y,z)\longmapsto\bigl(x,y,z-y\tilde\alpha(x,y)\bigr)
 \]
 gives a local diffeomorphism at the origin.
 Replacing $f$ by $\Phi_1\circ f$, we may set
 \[
    f(u,v)=\bigl(u,v^2,v^3 \tilde\beta(u,v^2)\bigr).
 \]
 Then by a straightforward calculation, 
 the unit normal vector field $\nu$ of $f$ is obtained as 
 \begin{multline*}
    \nu:=\frac{1}{\Delta}
      \left(-v^3\tilde\beta_u,
       -\frac32 v\tilde\beta-v^3 \tilde \beta_v,
      1\right),\\
   \Delta=
   \left[
      1+v^2
      \left(\left(\frac32 \tilde\beta+v^2 \tilde \beta_v \right)^2+v^{4}
      (\tilde\beta_u)^2 \right)\right]^{1/2}.
 \end{multline*}
 Since $\nu_v(u,0)=\bigl(0,-3\tilde\beta(u,0)/2,0\bigr)$, we have
 \[
    \tilde \psi(u)
            =\det(f_u,\nu_v,\nu)=
    \det
     \begin{pmatrix}
      1 & 0 & 0\\
      0 & -\frac{3}{2}\tilde\beta(u,0)& 0 \\
      0  & 0 & 1 
     \end{pmatrix}
       =-\frac{3}{2}\tilde\beta(u,0).
 \]
 Thus \ref{cond:main:2} of Theorem~\ref{main} holds if and only if
 \begin{equation}\label{eq:beta}
  \tilde\beta(0,0)=0,\qquad \tilde\beta_u(0,0)\ne 0.
 \end{equation}
 Then by the implicit function theorem,
 there exists a $C^\infty$-function
 $\delta(u,v)$ such that
 $\delta(0,0)=0$, and
 \begin{equation}\label{eq:delta}
   \tilde\beta\bigl(\delta(u,v),v\bigr)=u
 \end{equation}
 holds.
 Using this, we have a local diffeomorphism on $\R^2$ as
 $\phi:(u,v)\mapsto \bigl(\delta(u,v^2),v\bigr)$, 
 and
 \[
     f\circ \phi(u,v)=\bigl(\delta(u,v^2),v^2,uv^3\bigr).
 \]
 Since $\delta_u\ne 0$ by \eqref{eq:delta},
 $\Phi_2:(x,y,z)\mapsto (\delta(x,y),y,z)$
 gives a local diffeomorphism on $\R^3$, and
 \[
      \Phi_2^{-1}\circ f \circ \phi=(u,v^2,uv^3)
 \]
 gives the standard cuspidal cross cap $f_{\CCR}$ mentioned in 
the  introduction.
\end{proof}

In \cite{SUY}, the last three authors 
introduced the notion of {\em singular curvature\/}
of cuspidal edges, and studied the behavior of the Gaussian curvature
near a cuspidal edge:
\begin{fact}\label{fact:gauss}
 Let $f\colon{}U\to \R^3$ be a front,
 $p\in U$ a cuspidal edge,
 and $\gamma(t)$ $(|t|<\varepsilon)$ a singular curve consisting of
 non-degenerate singular  points with $\gamma(0)=p$.
 Then the Gaussian curvature $K$ is bounded on
 a sufficiently small neighborhood of
 $J:=\gamma\bigl((-\varepsilon,\varepsilon)\bigr)$
 if and only if the second fundamental form vanishes on 
 $J$.
 Moreover, 
 if the Gaussian curvature $K$
 is non-negative on $U\setminus J$ for a neighborhood of 
 $U$ of $p$, 
 then the singular curvature is non-positive.
\end{fact}
The singular curvature at a cuspidal cross cap is also 
defined in a similar way to the cuspidal edge case.
Since the unit normal vector field $\nu$ is well-defined at
a cuspidal cross caps, the second fundamental form is well-defined.
Since singular points sufficiently close to a cuspidal cross cap
are cuspidal edges, the following assertion immediately follows
from the above fact.
\begin{proposition}\label{prop:curvature}
 Let $f\colon{}U\to \R^3$ be a frontal,
 $p\in U$  a cuspidal cross cap,
 and $\gamma(t)$ $(|t|<\varepsilon)$ a singular curve consisting of
 non-degenerate singular  points with $\gamma(0)=p$.
 Then the Gaussian curvature $K$ is bounded on
 a sufficiently small neighborhood of
 $J:=\gamma\bigl((-\varepsilon,\varepsilon)\bigr)$
 if and only if the second fundamental form vanishes on 
 $J$.
 Moreover, 
 if the Gaussian curvature $K$
 is non-negative on $U\setminus J$ for a neighborhood of 
 $U$ of $p$, 
 then the singular curvature is non-positive.
\end{proposition}
Now, we give an example of a surface with umbilic points
accumulating at a cuspidal cross cap point.
For a space curve $\gamma(t)$ with arc-length parameter,
we take $\{\xi_1(t),\xi_2(t),\xi_3(t)\}$, 
$\kappa(t)>0$ and $\tau(t)$ as
the Frenet frame, the curvature and the torsion functions of $\gamma$.
We consider a tangent developable surface
$f(t,u)=\gamma(t)+u\xi_1(t)$ of $\gamma$.
The set of singular points of $f$ is $\{(t,0)\}$.

\begin{figure}
 \begin{center}
  \includegraphics[width=3cm]{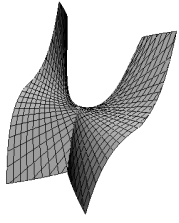}
 \end{center}
 \caption{Cuspidal cross cap with accumulating umbilic points}
 \label{fig:umbilic-line}
\end{figure}

We remark that this surface is frontal, since $\nu(t,u)=\xi_3(t)$
gives the unit normal vector.
By a direct calculation, the first fundamental form $ds^2$ and the
second fundamental form $\secondff$ are written as
\[
  ds^2 = \left(1+u^2\bigl(\kappa(t)\bigr)^2\right)\,dt^2+
                2\,dt\,du+du^2,\qquad
  \secondff = u\kappa(t)\tau(t)\,dt^2,
\] 
and the Gaussian curvature $K$ and the mean curvature are 
\[
     K=0, \qquad H= \frac{\tau(t)}{2u\kappa(t)}.
\]
So a regular point $(t,u)$ is an umbilic point if and only if 
$\tau(t)=0$.
On the other hand, it is easy to show that $f$ is a front at $(t,0)$
if and only if $\tau(t)\ne0$.
Moreover, Cleave~\cite{C} showed that a tangent developable surface $f$
at $(t,0)$ is ${\A}$-equivalent to a cuspidal cross cap if and only if
$\tau(t)=0$ and $\tau'(t)\ne 0$, which also follows from our 
criterion directly.
Hence we consider a tangent developable surface 
with space curve $\gamma(t)$ with $\tau(t)=0$ and $\tau'(t)\ne0$,
and then we have the desired example.
\begin{example}\label{ex:umbilic-line}
Let $\gamma(t)=(t,t^2,t^4)$ and consider a tangent developable
surface $f$ of $\gamma$.
Since $\tau(0)=0$ and $\tau'(0)=12\ne 0$,
all points on the ruling passing through 
$\gamma(0)$ are umbilic points
and $f$ at $(0,0)$ is a cuspidal cross cap 
(see Figure~\ref{fig:umbilic-line}). 
\end{example}

\section{Singularities of maximal surfaces}
\label{sec:maxface}
A holomorphic map $F=(F^1,F^2,F^3):M^2\to \C^3$ of  a Riemann surface 
$M^2$ into the complex space form $\C^3$ 
is called {\em null\/}
if $\sum_{j=1}^3 F^j_z\cdot F^j_z$ vanishes,
where $F^j_z:=dF^j/dz$.
We consider two projections, the former is the projection into 
the Euclidean $3$-space
\[
  p_E:\C^3\ni (\zeta^1,\zeta^2,\zeta^3)\mapsto 
      \Re(\zeta^1,\zeta^2,\zeta^3)\in\R^3,
\]
and the latter one is the projection into 
Lorentz-Minkowski $3$-space
\[
  p_L:\C^3\ni (\zeta^1,\zeta^2,\zeta^3)\mapsto 
      \Re(-\sqrt{-1} \zeta^3,\zeta^1,\zeta^2)\in L^3,
\]
where $L^3$ is the Lorentz-Minkowski space-time of dimension $3$
with signature $(-,+,+)$.
It is well-known that the projection of null holomorphic immersions
into $\R^3$ by $p_E$ gives conformal minimal immersions.
Moreover, conformal minimal immersions are always given locally
in such a manner.

On the other hand, the projection of null holomorphic immersions
into $L^3$ by $p_L$ gives spacelike maximal surfaces with singularities,
called {\it maxfaces} (see \cite{UY} for details).
Moreover, \cite{UY} proves that maxfaces are all frontal
and gives a necessary and sufficient condition for their 
singular points to be cuspidal edges and swallowtails.
In this section, we shall give a 
necessary and sufficient condition for their singular points
to be cuspidal cross caps and will show that generic singular points
of maxfaces consist of cuspidal edges, swallowtails and
cuspidal cross caps (see Theorem~\ref{thm:maxface} in 
the introduction).

The following fact is known (see \cite{UY}):
\begin{fact}
 Let $U\subset \C$ be a simply connected domain containing a base point
 $z_0$, and
 $(g,\omega)$ a pair of a meromorphic function and
 a holomorphic $1$-form on $U$ such that
 \begin{equation}\label{eq:normalized}
   (1+|g|^2)^2|\omega|^2\neq 0
 \end{equation}
 on $U$. Then 
 \begin{equation}\label{eq:maxface}
   f(z):=\Re\int_{z_0}^z \bigl(-2g, 1+g^2,\sqrt{-1}(1-g^2)\bigr)\omega
 \end{equation}
 gives a maxface in $L^3$.
 Moreover, any maxfaces are locally obtained in this manner.
\end{fact}

The first fundamental form (that is,  the induced metric)
of $f$ in \eqref{eq:maxface} is given by
\[
    ds^2=(1-|g|^2)^2|\omega|^2.
\]
In particular, $z\in U$ is a singular point of $f$ if and only if
$|g(z)|=1$, and 
at $f\colon{}U\setminus\{|g|=1\}\to L^3$
is a spacelike maximal (that is, vanishing mean curvature) immersion.
The meromorphic function $g$ can be identified with 
the Lorentzian Gauss map.
We call the pair $(g,\omega)$ the {\em Weierstrass data\/} of $f$.
In \cite{UY}, the last two authors proved that 
{\it $f$ is a front on a neighborhood of  a given singular point $z=p$ if
and only if\/} $\Re\bigl(dg/(g^2\omega)\bigr)\ne 0$. 
Moreover, the following assertions are proved in \cite{UY}: 

\begin{fact}[{\cite[Theorem 3.1]{UY}}]\label{fact:maxface-sing}
 Let $U$ be a domain of the complex plane $(\C,z)$ and
 $f:U \to L^3$ a maxface constructed from the Weierstrass data
 $(g,\omega=\hat\omega\, dz)$,
 where $\hat\omega$ is a holomorphic function on $U$. 
 Then $f$ is a frontal into $L^3$ {\rm (}which is identified with
 $\R^3${\rm )}.
 Take an arbitrary point $p\in U$.
 Then $p$ is a singular point of $f$ if and only if $|g(p)|=1$,
 and 
 $f$ is a front at a singular point $p$ if and only if
 $\Re\bigl(g'/(g^2\hat\omega)\bigr)\ne 0$ holds at $p$, 
 where $'=d/dz$. 
 Suppose now $\Re\bigl(g'/(g^2\hat\omega)\bigr)\ne 0$ at a singular
 point $p$.
 Then
  \begin{enumerate}
   \item $f$ is ${\A}$-equivalent to a cuspidal edge at $p$ if and only if
     $\Im\bigl(g'/(g^2\hat\omega)\bigr)\ne 0$, and
   \item $f$ is ${\A}$-equivalent to a swallowtail at $p$ if and only if
     \[
        \frac{g'}{g^2\hat\omega}\in\R\setminus\{0\}
        \qquad\text{and}\qquad
            \Re\left[
               \frac{g}{g'}\left(\frac{g'}{g^2\hat\omega}\right)'\right]
             \ne 0.
     \]
  \end{enumerate}
\end{fact}
\begin{remark}
 In \cite[Lemma~3.3]{UY},
 the third and fourth authors proved that 
 $\Re(g'/(g^2\hat\omega))\neq 0$
 at a singular point if and only if $f$ is a front at $p$ 
 and $p$ is {\em non-degenerate\/} singular point.
 More precisely, the condition $\Re(g'/(g^2\hat\omega))\neq 0$ 
 is equivalent to that $f$ is a front, and if this is the case,
 $p$ is automatically a non-degenerate singular point 
 because $g'\neq 0$.
 See page 25 in \cite{UY}.
\end{remark}
The statements of Theorem~3.1 of \cite{UY} are criteria 
to be {\em locally diffeomorphic\/} to a cuspidal edge
or a swallowtail.
However, in this case, local diffeomorphicity implies 
$\A$-equivalency.
See the appendix of \cite{KRSUY}.
We shall prove the following:

\begin{theorem}\label{thm:maxface-ccr}
 Let $U$ be a domain of the complex plane $(\C,z)$ and
 $f\colon{}U \to L^3$ a maxface constructed from the Weierstrass data
 $(g,\omega=\hat\omega\, dz)$,
 where $\hat\omega$ is a holomorphic function on $U$. 
 Take an arbitrary singular point  $p\in U$.
 Then 
 $f$ is ${\A}$-equivalent to a cuspidal cross cap at $p$
 if and only if 
 \[
    \frac{g'}{g^2\hat\omega}
     \in\sqrt{-1}\R\setminus\{0\}
    \qquad\text{and}\qquad
        \Im\left[
           \frac{g}{g'}\left(\frac{g'}{g^2\hat\omega}\right)'\right]
         \ne 0,
 \]
 where $'=d/dz$.
\end{theorem}
\begin{proof}
 We identify $L^3$ with the Euclidean $3$-space $\R^3$.
 Let $f$ be a maxface as in \eqref{eq:maxface}.
 Then 
 \[
    \nu:=\frac{1}{\sqrt{(1+|g|^2)^2+4 |g|^2}}
            (1+|g|^2,2\Re g,2\Im g)
 \]
 is the unit normal vector field of $f$ with respect to 
 the Euclidean metric of $\R^3$.
 Let $p\in U$ be a singular point of $f$, that is, 
 $|g(p)|=1$ holds.
 Then by \eqref{eq:normalized}, $\omega$ does not vanish 
 at $p$.
 Here, 
 \[
    \lambda = 
    \det (f_u,f_v,\nu) =
    (|g|^2-1)|\hat\omega|^2\sqrt{(1+|g|^2)^2+4|g|^2},
 \]
 under the complex coordinate $z=u+\sqrt{-1}v$ on $U$.
 Then the singular point $p$ is non-degenerate if and only if 
 $dg\neq 0$.

 The singular direction $\xi$ and the null direction $\eta$
 are given by
 $\xi=\sqrt{-1} \overline{(g'/g)}$, and
 $\eta = \sqrt{-1}/(g\hat\omega)$, respectively.
 Thus, we can parametrize the singular curve $\gamma(t)$ 
 as
 \begin{equation}\label{eq:dot-gamma}
    \dot\gamma(t) =\sqrt{-1}
   \overline{\left(\frac{g'}{g}\right)}\bigl(\gamma(t)\bigr)
   \qquad \left(\dot{~}=\frac{d}{dt}\right)
 \end{equation}
 (see the proof of Theorem~3.1 in \cite{UY}).
 Here, 
 we identify $T_pU$ with $\R^2$ and $\C$ with
 \begin{equation}\label{eq:equiv}
   \zeta=a+ \sqrt{-1} b \in \C
   ~\leftrightarrow~ (a,b)\in\R^2
   ~\leftrightarrow~ a\frac{\partial}{\partial u}+
           b\frac{\partial}{\partial  v}
   = \zeta\frac{\partial}{\partial z}+
           \bar \zeta\frac{\partial}{\partial \bar z},
 \end{equation}
 where $z=u+\sqrt{-1}v$.
 Then $\dot\gamma$ and $\eta$ are transversal if and only if 
 \[
   \det(\xi,\eta) = \Im \bar{\xi} \eta
                  = \Im \left(\frac{g'}{g^2\hat\omega}\right)\neq 0.
 \]

 On the other hand, one can compute $\psi$ as in
 Theorem~\ref{thm:criterion} as
 \[
   \psi = \det\bigl(f_*\dot\gamma,d\nu(\eta), \nu\bigr)
        = \Re\left(\frac{g'}{g^2\omega}\right)\cdot
          \psi_0,
 \]
 where $\psi_0$ is a smooth function on a neighborhood of $p$
 such that $\psi_0(p)\neq 0$.
 Then the second condition of Theorem~\ref{thm:criterion} is
 written as
 \[
   \Re\left(\frac{g'}{g^2\hat\omega}\right)=0\qquad\text{and}\qquad
   \Im\left[\left(\frac{g'}{g^2\hat\omega}\right)'
            \overline{\left(\frac{g'}{g}\right)}
           \right]\neq 0.
 \]
 Here, we used the relation 
 $d/dt = \sqrt{-1}[\overline{(g'/g)}(\partial/\partial z)
                         -(g'/g)(\partial/\partial\bar z)]$,
 which comes from \eqref{eq:dot-gamma}.
 Using  the relation
 $\overline{(g'/g)}  =(g/g')\cdot\mbox{(real valued function)}$,
 we have the conclusion.
\end{proof}
Let $f:U\to L^3$ be a maxface with a Weierstrass data $(g,\omega)$.
Then the associated maxface $\tilde f:U\to L^3$
with the Weierstrass data $(g,\sqrt{-1}\omega)$ is called
the {\it conjugate maxface}, which has the same first fundamental 
form and principal curvatures as $f$.
The following assertion follows immediately:
\begin{corollary}[%
 A duality between swallowtails and cuspidal cross caps]
 \label{cor:duality}
 A maxface $f:U\to L^3$ is $\A$-equivalent to a 
 swallowtail {\rm (}resp.\ a cuspidal cross cap{\rm)} at $p\in U$
 if and only if its conjugate maxface $\tilde f$
 is $\A$-equivalent to a cuspidal cross cap 
 {\rm (}resp.\ a swallowtail{\rm )} at $p\in U$.
\end{corollary}
\begin{remark*}
 The same assertion holds for \cmcone{} faces in de Sitter $3$-space,
 see Corollary~\ref{cor:desitter-duality}.
\end{remark*}
\begin{example}\label{ex:enneper}
 The Lorentzian Enneper surface is a maxface $f\colon{}\C\to L^3$
 with the Weierstrass data $(g,\omega)=(z,dz)$
 (see \cite[Example 5.2]{UY}), whose set of singularities
 is $\{z\,;\,|z|=1\}$.
 As pointed out in \cite{UY}, 
 Fact~\ref{fact:maxface-sing} implies that
 the points of the set 
 \[
    \{ z\,;\, |z|=1\} \setminus
    \bigl\{\pm 1, \pm \sqrt{-1}, \pm e^{\pm\sqrt{-1}\pi/4}\bigr\}
 \]
 are cuspidal edges and the points $\pm 1$, $\pm \sqrt{-1}$ are
 swallowtails.
 Moreover, using Theorem~\ref{thm:maxface-ccr},
 we deduce that the four points $\pm e^{\pm \sqrt{-1}\pi/4}$ are
 cuspidal cross caps.
\end{example}
By \eqref{eq:normalized}, $\omega$ does not vanish at 
a singular point $p$.
Hence there exists a complex coordinate system
$z$ such that $\omega=dz$.
On the other hand, $g\neq 0$ at the singular point $p$.
Hence there exists a holomorphic function $h$ in $z$ 
such that $g=e^h$.
We denote by $f_h$ the maxface defined by 
the Weierstrass data $(g,\omega)=(e^h,dz)$.
Let $\O(U)$ be the set of holomorphic functions defined on $U$,
which is endowed with the compact open $C^{\infty}$-topology.
Then we have the induced topology on  the set of maxfaces 
$\{f_h\}_{h\in \O(U)}$.
We shall prove Theorem~\ref{thm:maxface} in the introduction.
To prove the theorem, we rewrite the criteria in
Fact~\ref{fact:maxface-sing} and Theorem~\ref{thm:maxface-ccr}
in terms of $h$.
\begin{lemma}\label{lem:maxface-criteria}
 Let $h\in \O(U)$ and set
 \[
     \alpha_h:= e^{-h}h',\qquad
     \beta_h:= e^{-2h}\bigl(h''-(h')^2\bigr),
 \]
 where $'=d/dz$.
 Then 
 \begin{enumerate}
  \item\label{maxface:1}
       a point $p\in U$ is a singular point of $f_h$ if and only
       if  $\Re h=0$,
  \item\label{maxface:2}
        a singular point $p$ is non-degenerate if and only if
    $\alpha_h\neq 0$, 
  \item\label{maxface:3}
        a singular point $p$ is a cuspidal edge
    if and only if $\Re \alpha_h\neq 0$ and $\Im \alpha_h\neq 0$,
  \item\label{maxface:4} 
        a singular point $p$ is a swallowtail
    if and only if $\Re\alpha_h\neq 0$,
    $\Im \alpha_h=0$, and $\Re \beta_h\neq 0$,
  \item\label{maxface:5} 
        a singular point $p$ is a cuspidal cross cap
    if and only if $\Re\alpha_h= 0$,
    $\Im \alpha_h\neq 0$, and $\Re \beta_h\neq 0$.
 \end{enumerate}
\end{lemma}
\begin{proof}
 Since $g=e^h$, \ref{maxface:1} is obvious.
 Moreover, a singular point $p$ is non-degenerate
 if and only if $g'=e^{2h}\alpha_h$ does not vanish.
 Hence we have \ref{maxface:2}.
 Since $g'/(g^2\hat\omega) = \alpha_h$, the criterion for 
 a front (Fact~\ref{fact:maxface-sing}) 
 is $\Re \alpha_h\neq 0$.
 Then by Fact~\ref{fact:maxface-sing}, 
 we have \ref{maxface:3}.
 Here,
 \[
     \frac{g}{g'}\left(\frac{g'}{g^2\hat\omega}\right)'
     = \frac{e^{-h}}{h'}\left(h''-(h')^2\right).
 \]
 Then, if $\Im\alpha_h=0$ and $\alpha_h\neq 0$, 
 \[
    \Re\left[\frac{g}{g'}\left(\frac{g'}{g^2\hat\omega}\right)'\right]
    =
    \Re\left[
        e^{-2h}\left(h''-(h')^2\right)
       \right]\frac{1}{\alpha_h}
    =\frac{1}{\alpha_h}\Re\beta_h.
 \]
 Then by Fact~\ref{fact:maxface-sing}, we have \ref{maxface:4}.
 On the other hand, if $\Re\alpha_h=0$ and $\alpha_h\neq 0$,
 \[
    \Im\left[\frac{g}{g'}\left(\frac{g'}{g^2\hat\omega}\right)'\right]
    =
    \Im\left[
      \frac{e^{-2h}}{\alpha_h}\left(h''-(h')^2\right)\right]
      =-\frac{\sqrt{-1}}{\alpha_h}\Re\beta_h.
 \]
Thus we have \ref{maxface:5}.
\end{proof}
Let $J^2_{\O}(U)$ be the space of $2$-jets of holomorphic functions 
on $U$, which is identified with an $8$-dimensional manifold
\[
  J^2_{\O}(U)=U\times \C\times \C\times \C
          =U\times \F \times \F_1 \times \F_2,
\]
where $\F$, $\F_1$ and $\F_2$ 
are $1$-dimensional complex vector spaces corresponding to 
the jets 
$h(p)$, $h'(p)$ and $h''(p)$ for $h\in \O(U)$ at $p$, respectively.
Then, the canonical map $j^2_h\colon{}U\to J^2_{\O}(U)$
is given by $j^2_h(p) = \bigl(p,h(p),h'(p),h''(p) \bigr)$.
The point $P\in J^2_{\O}(U)$ is expressed as 
\begin{equation}\label{eq:holo-parameter}
  P=(p,\hat h,\hat h_1,\hat h_2)
   =(p, \hat u, \hat v, \hat u_1, \hat v_1, \hat u_2, \hat v_2),
\end{equation}
where $\hat h = \hat u+\sqrt{-1}\hat v$, 
      $\hat h_j=\hat u_j+\sqrt{-1}\hat v_j$ ($j=1,2$).
We set 
\begin{align*}
   A&:=\{P\in J^2_{\O}(U)\,;\, \Re \hat h=0, \hat \alpha=0\},\\
   B&:=\{P\in J^2_{\O}(U)\,;\, \Re \hat h=0, \Im \hat\alpha=0,
                            \Re\hat\beta=0\},\\
   C&:=\{P\in J^2_{\O}(U)\,;\, \Re \hat h=0, \Re \hat\alpha=0,
                            \Re\hat\beta=0\},
\end{align*}
where 
\[
  \hat\alpha = e^{-\hat h}\hat h_1,\qquad
  \hat\beta = e^{-2\hat h}(\hat h_2-(\hat h_1)^2).
\]
\begin{lemma}\label{lem:criteria-h}
 Let $S=A\cup B\cup C$ and 
 \[
    \G := \{ h\in \O(U)\,;\,j^{2}_h(U)\cap S=\emptyset \}.
 \]
 Then all singular points of $f_h$ are cuspidal edges,
 swallowtails or cuspidal cross caps if $h\in\G$.
\end{lemma}
\begin{proof}
 We set 
 \begin{align*}
   \S_A&:=\{h\in \O(U)\,;\,j^{2}_h(U)\cap A\neq\emptyset\},\\
   \S_B&:=\{h\in \O(U)\,;\,j^{2}_h(U)\cap B\neq\emptyset\},\\
   \S_C&:=\{h\in \O(U)\,;\,j^{2}_h(U)\cap C\neq\emptyset\}.
 \end{align*}
 Then we have
 $\G=(\S_A)^c\cap(\S_B)^c\cap(\S_C)^c$.
 Let $h\in \G$, and let $p\in U$ be a singular point of $f_h$.
 Since $h\not\in\S_A$, $p$ is a non-degenerate singular point.
 If $f_h$ is not a front at $p$, then $\Re\alpha_h=0$.
 Since $h\not\in\S_C$, this implies that $\Re\beta_h\neq 0$,
 and hence $p$ is a cuspidal cross cap.
 If $f_h$ is a front at $p$ and not a cuspidal edge,
 $p$ is a swallowtail since $h\not\in\S_B$.
\end{proof}

By this lemma, a singular point $p\in U$ is 
neither a cuspidal edge,
a swallowtail nor a cuspidal cross cap if and only if $j^2_h(p)\in S$.
Thus if $S$ is a union of a finite number of
 submanifolds in $J^2_{\O}(U)$ of codimension $3$,
then we can say generic singular points of $h$ 
consist of cuspidal edges, swallowtails and 
cuspidal cross caps.  
In fact, Theorem~\ref{thm:maxface} can be proved 
in a similar way to Theorem~3.4 of \cite{KRSUY}
using the following lemma.

\begin{lemma}\label{lem:codim-maxface}
 $S=A\cup B\cup C$ is the  union of a finite number of
 submanifolds in $J^2_{\O}(U)$ of codimension $3$.
\end{lemma}
\begin{proof}
 Using parameters in \eqref{eq:holo-parameter},
 we can write
 \[
     A = \{ \hat u=\hat u_1=\hat v_1=0\},
 \]
 which is a codimension $3$ submanifold in $J^2_{\O}(U)$.
 Moreover, one can write 
 \begin{align*}
    B&=\{ \zeta_1 = 0 , \zeta_2 = 0 , \zeta_3 = 0 \},
    \quad\text{where}\\
        &\zeta_1= \hat u,\\
        &\zeta_2= e^{-\hat  u}(\hat v_1\cos \hat v-\hat u_1\sin\hat v),\\
        &\zeta_3= e^{-2\hat u}
             \biggl(\left(\hat u_2-(\hat u_1)^2+(\hat v_1)^2\right)\cos 2\hat v+
                   (\hat v_2-2\hat u_1\hat v_1)\sin 2\hat v\biggr).
 \end{align*}
 Here, we can compute that
 \begin{equation}\label{eq:jacobi}
   \frac{%
       \partial (\zeta_1,\zeta_2,\zeta_3)}{%
       \partial (\hat u,\hat u_1,\hat v_1)} = 
        2 e^{-3\hat u}(\hat u_1 \cos\hat v + \hat v_1 \sin\hat v).
 \end{equation}
 Since $(\hat u_1,\hat v_1)\neq (0,0)$ and  
 $\hat v_1 \cos \hat v-\hat u_1 \sin \hat v=0$ hold
 on $B\setminus A$,  
 \eqref{eq:jacobi} does not vanish on $B\setminus A$.
 Hence by the implicit function theorem, $B\setminus A$ is a submanifold
 of codimension $3$.

  Similarly, $C$ is written as
  \begin{align*}
     C &= \{\xi_1=0,\xi_2=0,\xi_3=0\},    \quad\text{where}\\
       &\xi_1 = \hat u,\qquad \xi_2 = 
           e^{-\hat u}(\hat u_1\cos \hat v+ \hat v_1 \sin \hat v),\\
       &\xi_3= e^{-2\hat u}
             \biggl(\left(\hat u_2-(\hat u_1)^2+(\hat v_1)^2\right)\cos 2\hat v+
                   (\hat v_2-2\hat u_1\hat v_1)\sin 2\hat v \biggr).
  \end{align*}
  Then we have 
  \[
       \frac{\partial(\xi_1,\xi_2,\xi_3)}{%
       \partial(\hat u,\hat u_1,\hat v_1)}
         =
         2 e^{-3\hat u}(\hat v_1 \cos\hat v - \hat u_1 \sin\hat v).
  \]
  Thus, $C\setminus A$ is a submanifold of codimension $3$.

  Hence $S = A \cup (B\setminus A)\cup(C\setminus A)$ is a 
  union of submanifolds of codimension $3$.
\end{proof}

\section{Singularities of \cmcone{} surfaces in de Sitter 3-space}
\label{sec:cmc1-face}

It is well-known that spacelike
\cmcone{} (constant mean curvature one)
surfaces in de Sitter 3-space $S^3_1$
have similar properties to spacelike maximal surfaces in $L^3$.
In this section, we shall give an analogue of 
Theorem~\ref{thm:maxface} for such surfaces.  
Though the assertion is the same, the method is not parallel:
For maxfaces, one can easily 
write down the Euclidean normal vector explicitly, 
as well as the Lorentzian normal, in terms of the Weierstrass data.
However, the case of \cmcone{} surfaces in $S^3_1$ is different,
as it is difficult to express the Euclidean normal vector, and we 
apply Theorem~\ref{main} instead of Corollary~\ref{thm:criterion}, 
since Theorem~\ref{main} is independent of the metric of the ambient space.

A holomorphic map $F\colon{}M^2\to \SL(2,\C)$ of  a Riemann surface 
$M^2$ into the complex Lie group $\SL(2,\C)$ 
is called {\em null\/}
if $\det F_z=0$ holds on $M^2$, where $z$ is a local complex coordinate 
of $M^2$.
This condition does not depend on the choice of complex coordinates.
We consider two projections of holomorphic null immersions, the 
former is
the projection into the hyperbolic $3$-space
\[
  p_E:\SL(2,\C)\ni F \mapsto FF^*\in H^3=\SL(2,\C)/\SU(2)
\]
and the latter one is the projection into 
de Sitter $3$-space $S^3_1$
\[
  p_L:\SL(2,\C)\ni F \mapsto 
  F\begin{pmatrix} 1 & \hphantom{-}0 \\ 0 & -1 \end{pmatrix}F^*
  \in S^3_1=\SL(2,\C)/\SU(1,1),
\]
where  $F^*={}^t\overline F$.
It is well-known that the projection of null holomorphic immersions
into $H^3$ by $p_E$ gives conformal \cmcone{} immersions
(see \cite{Br} and \cite{UY1}).
Moreover, conformal \cmcone{} immersions are always given locally
in such a manner.

On the other hand, the projection of null holomorphic immersions
into $S^3_1$ by $p_L$ gives spacelike \cmcone{} surfaces with
singularities, called {\it \cmcone{} faces} (see \cite{F} for details).
In this section, we shall give a necessary and sufficient condition for
their singular points
to be cuspidal edges, swallowtails and cuspidal cross caps, and will show
that \cmcone{} faces admitting only  these singular points are generic.

Recall that de Sitter $3$-space is 
\[
  S^3_1=S^3_1(1)=
  \{(x^{0},x^{1},x^{2},x^{3})\in L^4 \, ; \, 
   -(x^{0})^2+(x^{1})^2+(x^{2})^2+(x^{3})^2=1\},
\]
with metric induced from $L^4$, which is 
a simply-connected $3$-dimensional Lorentzian manifold with constant
sectional curvature $1$. 
We can consider $L^4$ to be the set of $2 \times 2$ 
Hermitian matrices $\Herm(2)$ by the identification 
\[
   L^4\ni X=(x^{0},x^{1},x^{2},x^{3})\leftrightarrow 
   X=\sum_{k=0}^3 x^{k} e_k 
   =\begin{pmatrix} x^{0}+x^{3} & 
                  x^{1}+\sqrt{-1} x^{2} \\ 
                  x^{1}-\sqrt{-1} x^{2} & 
                  x^{0}-x^{3} 
    \end{pmatrix},
\]
where
\[
      e_0=\begin{pmatrix}1&0\\0&1\end{pmatrix},\quad
      e_1=\begin{pmatrix}0&1\\1&0\end{pmatrix},\quad
      e_2=\begin{pmatrix}0&\sqrt{-1}\\-\sqrt{-1}&0\end{pmatrix},\quad
      e_3=\begin{pmatrix}1&0\\0&-1\end{pmatrix}.
\]
Then $S^3_1$ is written as
\[ 
   S^3_1=\{X\,;\,X^*=X\,,\det X=-1\}
                  =\{Fe_3F^*\,;\,F\in\SL(2,\C)\} 
\] 
with the metric
\[
 \langle X,Y\rangle 
      =-\frac{1}{2}\trace\left(Xe_2({}^t\!{Y})e_2\right) . 
\]
In particular, $\langle X,X\rangle =-\det X$. 

The following fact is known (see \cite{F}):
\begin{fact}
 Let $U\subset \C$ be a simply connected domain containing a base point
 $z_0$, and
 $(g,\omega)$ a pair of a meromorphic function and
 a holomorphic $1$-form on $U$ such that 
\eqref{eq:normalized} holds on $U$. 
 Then by solving the ordinary differential equation
 \begin{equation}\label{eq:F^-1dF}
  F^{-1}dF=
   \begin{pmatrix}
    g & -g^2 \\ 1 & -g\hphantom{^2}
   \end{pmatrix} \omega
 \end{equation}
 with $F(z_0)=e_0$, 
 where $z_0\in U$ is the base point,
 \begin{equation}\label{eq:cmc1-face}
   f:=Fe_3F^*
 \end{equation}
 gives a \cmcone{} face in $S^3_1$.
 Moreover, any \cmcone{} face is locally obtained in this manner.
\end{fact}
The first fundamental form 
of $f$ in \eqref{eq:cmc1-face} is given by
\[
    ds^2=(1-|g|^2)^2|\omega|^2.
\]
In particular, $z\in U$ is a singular point of $f$ if and only if
$|g(z)|=1$, and 
 $f\colon{}U\setminus\{|g|=1\}\to S^3_1$
is a spacelike \cmcone{} (constant mean curvature one) immersion.
We call the pair $(g,\omega)$ the {\em Weierstrass data\/} of $f$ and
$F$ the {\it holomorphic null lift} of $f$.

Let $f:U\to S^3_1$ be a \cmcone{} face and $F$ a holomorphic
null lift of $f$ with Weierstrass data $(g,\omega)$.
We set
\[
  \beta:=
  \begin{pmatrix}
   1 & g \\ \bar g & 1
  \end{pmatrix}.
\]
Then
\begin{equation}\label{eq:cmcone-nu}
   \nu:=F \beta^2 F^*  
\end{equation}
gives the Lorentzian normal vector field of $f$ on 
the regular set of $f$.
This $\nu$ is not a unit vector, 
but extends smoothly across the singular sets.
Let $TS^3_1|_{f(U)}$ be the restriction of the tangent bundle of $S^3_1$
to $f(U)$. Then 
\[
  L:=\langle {*},{\nu}\rangle
\]
gives a section of $U$ into $T^*U$, and gives an admissible lift of $f$.
In particular, $f$ is a frontal and
the subbundle
\[
   \{X\in TS^3_1|_{f(U)}\,;\, \langle {X},{\nu}\rangle =0 \}
\]
coincides with the limiting tangent bundle.
Moreover, we have the following:

\begin{lemma}
 Any section $X$ of the limiting tangent bundle is parametrized
 as 
 \begin{equation}\label{eq:paramE}
  X=F\begin{pmatrix}\bar\zeta g+\zeta\bar g&\zeta (|g|^2+1)\\
     \bar\zeta (|g|^2+1)&\bar\zeta g+\zeta\bar g\end{pmatrix}F^*
 \end{equation}
 for some $\zeta :U\to\C$. 
\end{lemma}
\begin{proof}
 Let $p$ be an arbitrary point in $U$. 
 Since $X_p\in \Herm(2)$, 
 $X_p\in T_p L^4$. 
 Because $\inner{f_p}{X_p}=0$, $X_p\in T_pS^3_1$. 
 Since $\inner{\nu_p}{X_p}=0$, 
 and $\inner{~}{~}$ is a non-degenerate inner product, 
 we get the conclusion. 
\end{proof}

The above lemma will play a crucial role in
giving a criterion for 
cuspidal cross caps in terms of the Weierstrass data.

\begin{proposition}\label{prop:non-deg}
Let $U$ be a domain of the complex plane $(\C,z)$. 
Let $f:U\to S^3_1$ be a \cmcone{} face and $F$ a holomorphic null lift 
of $f$ with Weierstrass data  $(g,\omega=\hat\omega\,dz)$, 
where $\hat\omega$ is a holomorphic function on $U$.  
Then a singular point $p\in U$ is 
non-degenerate if and only if $dg(p)\ne 0$. 
\end{proposition}

\begin{proof}
 Define $\xi\in T_{f(p)} L^4$ as 
\begin{equation}\label{eq:nu}
   \xi:=FF^*.
\end{equation}
 Then $\xi\in T_{f(p)}S^3_1$,  because $\inner{f}{\xi}=0$. 
 Define a $3$-form $\Omega$ on $S^3_1$ as 
 \begin{equation}\label{eq:Omega}
  \Omega (X_1,X_2,X_3):=\det (f,X_1,X_2,X_3)
 \end{equation} 
 for arbitrary vector fields $X_1,X_2,X_3$ of $S^3_1$. 
 Then $\Omega$ gives a volume element on $S^3_1$. 
 Since 
 \begin{align*}
    \Omega (f_{u},f_{v},\xi)
  &=\det (f,f_{u},f_{v},\xi) \\
  &=\det\begin{pmatrix}
    0&2\Re(g\hat\omega)        &-2\Im(g\hat\omega)       &1\\
    0& \Re\{(1+g^2)\hat\omega\}&-\Im\{(1+g^2)\hat\omega\}&0\\
    0&-\Im\{(1-g^2)\hat\omega\}&-\Re\{(1-g^2)\hat\omega\}&0\\
    1&  0                      &  0                      &0\end{pmatrix} \\
  &=(1-|g|^2)(1+|g|^2)|\hat\omega|^{2}, 
 \end{align*}
 we see that 
 \begin{align*}
  d\bigl(\Omega (f_{u},f_{v},\nu)\bigr)
  &=-\frac{1}{2}
    \biggl(d(g\bar g)(1+|g|^2)|\hat\omega|^{2}
                  -(1-|g|^2)d\bigl((1+|g|^2)|\hat\omega|^{2}\bigr)
    \biggr)\\
  &=-d(g\bar g)|\hat\omega|^{2}
 \end{align*}
 at $p$, because $|g(p)|=1$, proving the proposition 
 by Proposition~\ref{non-deg}.
\end{proof}

We shall now prove the following:

\begin{theorem}\label{thm:sing}
 Let $U$ be a domain of the complex plane $(\C,z)$ and
 $f:U \to S^3_1$ a \cmcone{} face constructed from the Weierstrass data
 $(g,\omega=\hat\omega\, dz)$,
 where $\hat\omega$ is a holomorphic function on $U$.
 Then:
 \begin{enumerate}
   \item\label{item:sing-1}
     A point $p\in U$ is a singular
     point if and only if $|g(p)|=1$.
   \item\label{item:sing-2}
    $f$ is ${\A}$-equivalent to a cuspidal edge at a singular point $p$
    if and only if
     \[
         \Re\left(\frac{g'}{g^2\hat\omega}\right)\neq 0
         \qquad\text{and}\qquad
         \Im\left(\frac{g'}{g^2\hat\omega}\right)\neq 0
     \]
     hold at $p$, where $'=d/dz$.
   \item\label{item:sing-3}
    $f$ is ${\A}$-equivalent to a swallowtail at a singular point $p$
    if and only if
     \[
        \frac{g'}{g^2\hat\omega}\in\R\setminus\{0\}\qquad
        \text{and}\qquad
        \Re\left\{
        \frac{g}{g'}\left(
           \frac{g'}{g^2\hat\omega}
        \right)'\right\}\neq 0
     \]
     hold at $p$.
   \item\label{item:sing-4}
    $f$ is ${\A}$-equivalent to a cuspidal cross cap at a singular point $p$
    if and only if
     \[
        \frac{g'}{g^2\hat\omega}\in\sqrt{-1}\R\setminus\{0\}
        \qquad\text{and}\qquad
        \Im\left\{
        \frac{g}{g'}\left(
           \frac{g'}{g^2\hat\omega}
        \right)'\right\}\neq 0
     \]
     hold at $p$.
 \end{enumerate}
In particular,
 the criteria for cuspidal edges, swallowtails and cuspidal cross 
 caps in terms of $(g,\omega)$ are exactly the same as in
 the case of maxfaces
 {\rm (}Fact~\ref{fact:maxface-sing} and Theorem~\ref{thm:maxface-ccr}{\rm)}.
\end{theorem}

To prove Theorem \ref{thm:sing}, 
we prepare the following lemmas:
\begin{lemma}\label{lem:front-cmc}
 Let $f:M^2\to S^3_1$ be a \cmcone{} face constructed from the
 Weierstrass data  $(g,\omega)$.
 Assume that $p\in M^2$ is a singular point. 
 Then $f$ is a front on a neighborhood of $p$ if and only if 
 \begin{equation}\label{eq:Re(dg/g^2w)}
  \Re\left(\frac{dg}{g^2\omega}\right) \neq 0 
 \end{equation}
 holds at $p$.  
 If this is the case, $p$ is a non-degenerate singular point.
\end{lemma}
\begin{proof}
 Let $\nu$ be as in \eqref{eq:cmcone-nu}.
 Then by the similar argument in the proof of Lemma~\ref{lem:front-2},
 $f$ is a front at $p$ if and only if $D_{\eta}^{L^4}\nu$ 
 is not proportional to $\nu$, where $\eta$ is the null direction 
 which is given by
\begin{equation}\label{eq:null}
 \eta=
  \frac{\sqrt{-1}}{g\hat \omega}\frac{\partial}{\partial z}
  -\frac{\sqrt{-1}}{\bar g \overline{\hat \omega}} 
  \frac{\partial}{\partial \bar z} 
\end{equation}
 at a singular point $p$, where $z$ is a complex coordinate and 
 $\omega=\hat\omega\,dz$.
 Using \eqref{eq:F^-1dF} and noticing $|g|=1$ at $p$, we have
 \begin{align*}
    D_\eta^{L^4}\nu &=
         \frac{\sqrt{-1}}{g\hat\omega}\frac{\partial}{\partial z}
         (F\beta^2F^* )
         -\frac{\sqrt{-1}}{\bar g \overline{\hat \omega}} 
         \frac{\partial}{\partial \bar z} 
         (F\beta^2F^*) \\
         &= 
          \frac{\sqrt{-1}}{g\hat\omega} 
           F\bigl(F^{-1}F_z\beta^2 + (\beta^2)_z\bigr)F^*-
          \frac{\sqrt{-1}}{\bar g \bar{\hat\omega}} 
           F\bigl(\beta^2(F^*)_{\bar z}(F^{*})^{-1} + (\beta^2)_{\bar z}\bigr) F^*
         \\
         &=
          2F\begin{pmatrix}
		         \Im g\bar\mu & \sqrt{-1}\mu \\
			  -{\sqrt{-1}}\bar\mu   & \Im g\bar \mu
	    \end{pmatrix}F^*,\qquad
                      \text{where}\quad
                      \mu = \frac{g'}{g\hat\omega}
 \end{align*}
 at $p$, where $D^{L^4}$ is the canonical connection of $L^4$.
 On the other hand, we have
 \[
   \inner{D_\eta^{L^4}{\nu}}{\nu} = 0.
 \]
 Here, $\inner{\nu}{\nu}=0$ at $p$ (that is, $\nu$ is a null vector
 in $L^4$).
 Then  $D_\eta^{L^4}\nu$ is proportional to $\nu$ if and only if 
 it is a null vector, which is equivalent to
 \begin{align*}
   0=\det D_\eta^{L^4}\nu &=
     4 \bigl((\Im g\bar \mu)^2    -\mu\bar \mu\bigr)
     =
     4 \bigl((\Im g\bar \mu)^2    -g\bar g\mu\bar \mu\bigr)\\
     &= -4(\Re g\bar \mu)^2=-4\left(\Re \bar g\mu\right)^2
     =-4\left(\Re \frac{\mu}{g}\right)^2 =
      -4\left(\Re \frac{g'}{g^2\hat\omega}\right)^2,
\end{align*}
 because $g\bar g=1$.
 Thus we have the conclusion.
 Moreover, if this is the case, $dg\neq 0$ at $p$.
 Then by Proposition~\ref{prop:non-deg}, $p$ is a non-degenerate
 singular point.
\end{proof}

\begin{lemma}\label{lm:ccr}
 Let $U$ be a domain of the complex plane $(\C,z)$.
 Let $f:U\to S^3_1$ be a \cmcone{} face 
 and $F$ a holomorphic null lift of $f$ 
 with Weierstrass data  $(g,\omega=\hat\omega)$, 
 where $\hat\omega$ is a holomorphic function on $U$.
 Let $X$ be a section of the limiting tangent bundle
 defined as in equation \eqref{eq:paramE}. 
 Take a singular point $p\in U$. 
 Then 
\[
 \tilde\psi:=
 \left\langle\nu,D^{f}_\eta X\right\rangle
 =2\Re\left(\frac{g'}{g^2\hat\omega}\right)\Im(\bar\zeta g)
\]
 holds, where $D$ is the canonical connection of $S^3_1$,
$\nu$ is the vector field given in \eqref{eq:cmcone-nu},
 and $\eta$ denotes the null direction of $f$ at $p$. 
\end{lemma}
\begin{proof}
 We set
 \[
 T=\begin{pmatrix}
    \bar\zeta g+\zeta\bar g&\zeta (|g|^2+1)\\
    \bar\zeta (|g|^2+1)&\bar\zeta g+\zeta\bar g
   \end{pmatrix}. 
\]
Then $X=FTF^*$. 
On the other hand, the null direction $\eta$ is given as
\eqref{eq:null} at $p$.
Thus 
\[
   D^{L^4}_\eta X\\
  =
     \frac{\sqrt{-1}}{g\hat \omega}F
      \bigl(F^{-1}F_zT + T_z\bigr)F
     -\frac{\sqrt{-1}}{\bar g \overline{\hat \omega}}
      F\bigl(F^{-1}F_zT)^* + T_{\bar z} \bigr) F^*.
\]
Since $\bar g=g^{-1}$ at any singular point $p$, 
and by \eqref{eq:F^-1dF}, we see that 
\begin{align*}
 \frac{\sqrt{-1}}{g\hat \omega}F(F^{-1}F_zT)F^*
 &=-\frac{\sqrt{-1}}{\bar g \overline{\hat \omega}}
 F(F^{-1}F_zT)^*F^* 
 =\sqrt{-1}F
 \begin{pmatrix}
   \zeta \bar g -\bar \zeta g & \zeta -\bar \zeta g^2 
  \\
  \zeta \bar g^2 - \bar \zeta & \zeta \bar g - \bar\zeta g 
 \end{pmatrix}F^* \\
 &=\sqrt{-1}(\zeta \bar g-\bar \zeta g)
  F\begin{pmatrix}1 & g \\ \bar g & 1\end{pmatrix}F^*. 
\end{align*}
Thus 
\[
  \inner{D^{L^4}_\eta X}{\nu}
    =\inner{\frac{\sqrt{-1}}{g\hat \omega}FT_zF^*
    -\frac{\sqrt{-1}}{\bar g \overline{\hat \omega}}FT_{\bar z}
        F^*}{\nu}.
\]
Since 
\begin{align*}
 \inner{\frac{\sqrt{-1}}{g\hat \omega}FT_zF^*}{\nu}
 &=-\frac{1}{2}\trace\left[\frac{\sqrt{-1}g'}{g\hat\omega}
   \begin{pmatrix}
    \bar\zeta&\zeta\bar g\\
   \bar\zeta\bar g&\bar\zeta
   \end{pmatrix}
   \begin{pmatrix}1&-g\\-\bar g&1\end{pmatrix}\right] \\
 &= \frac{g'}{2g^2\hat\omega}\sqrt{-1}(\zeta\bar g-\bar\zeta g) 
\end{align*}
and 
\begin{align*}
 \inner{\frac{\sqrt{-1}}{\bar g\overline{\hat\omega}}
 FT_{\bar z}F^*}{\nu}
 &=-\frac{1}{2}\mathrm{trace}
   \left[\frac{\sqrt{-1}g'}{\bar g\overline{\hat\omega}}
   \begin{pmatrix}\zeta&\zeta g\\\bar\zeta g&\zeta\end{pmatrix}
   \begin{pmatrix}1&-g\\-\bar g&1\end{pmatrix}\right] \\
 &=-\overline{\left(\frac{g'}{2g^2\hat\omega}\right)}
   \sqrt{-1}(\zeta\bar g-\bar\zeta g), 
\end{align*}
we have 
\[
 \tilde\psi
   =\inner{D^{L^4}_\eta X}{\nu}
    =2\Re\left(\frac{g'}{g^2\hat\omega}\right)
      \Im(\bar\zeta g), 
\]
proving the lemma. 
\end{proof}

Now assume that $X$ defined as in \eqref{eq:paramE} satisfies 
\ref{ass:main:2} in Theorem~\ref{main}. 
Then by the definition of $X$, $\Im(\bar\zeta g)$ cannot be zero at a 
singular point.  
\begin{proof}[Proof of Theorem~\ref{thm:sing}]
 Since the criteria for cuspidal edges and swallowtails are
 described intrinsically, and the first fundamental form
 of $f$ is the same as in the case of maxfaces,
 so the assertions 
\ref{item:sing-1}, \ref{item:sing-2} and \ref{item:sing-3} 
 are parallel to the case of maxfaces in $L^3$.  
 So it is sufficient to show the last assertion:
 Let $\gamma$ be the singular curve with $\gamma(0)=p$.
 Since the induced metric $ds^2$ is in the same form as for the
 maxface case,
 we can parametrize $\gamma$ as \eqref{eq:dot-gamma}:
\[
    \dot\gamma(t) = \sqrt{-1}
    \overline{
    \left(\frac{g'}{g}\right)}\bigl(\gamma(t)\bigr),
\]
where $\dot{~}=d/dt$.

On the other hand, the null direction is given as in \eqref{eq:null}.
Assume $X$ satisfies \ref{ass:main:2} in Theorem~\ref{main}. 
Then the necessary and sufficient condition for a cuspidal cross cap is
$\tilde\psi=0$ and $d\tilde\psi/dt \ne 0$, by Theorem~\ref{main}.
  Thus, Lemma~\ref{lm:ccr} implies the last assertion, since
  \begin{align*}
      \left.\frac{d}{dt} \right|_{t=0}\!\!  \inner{D^{L^4}_\eta X}{\nu} 
      &= 2\Im(\bar\zeta g)\Re\left[\left(\frac{g'}{g^2\hat\omega}\right)'
          \frac{d\gamma}{dt}\right] \\
      &=-2\Im(\bar\zeta g)\Im\left[\left(\frac{g'}{g^2\hat\omega}\right)'
            \overline{\left(\frac{g'}{g}\right)}\right] \\
      &=-2|g'|^2\Im(\bar\zeta g)
        \Im\left[\left(\frac{g'}{g^2\hat\omega}\right)'
        \frac{g}{g'}\right]. \qquad \qed
  \end{align*}
\renewcommand{\qedsymbol}{\relax}
\end{proof}

We take a holomorphic function $h$ defined on a simply connected
domain $U\subset\C$.
Then there is a \cmcone{} face $f_h$ with Weierstrass data 
$(g,\omega)=(e^h,dz)$,
where $z$ is a complex coordinate of $U$.
Let $\O(U)$ be the set of holomorphic functions on $U$,
which is endowed with the compact open $C^{\infty}$-topology.
Since the
criteria for cuspidal edges, swallowtails and cuspidal cross 
caps in terms of $(g,\omega)$ are exactly the same as 
in the case of maxfaces,
we have the following:
\begin{corollary}\label{thm:cmc1-face}
 Let $U\subset\C$ be a simply connected domain and $K$ an 
 arbitrary compact set,
 and let  $S(K)$ be the subset of  $\O(U)$ consisting of 
 $h\in \O(U)$
 such that the singular points of $f_h$ are cuspidal edges,
 swallowtails or cuspidal cross caps.
 Then $S(K)$ is an open and dense subset  of $\O(U)$.
\end{corollary}
As the same as case of maxfaces,
the {\em conjugate\/} \cmcone{} face $\tilde f$
of a \cmcone{} face $f$ is defined by the Weierstrass data
$(g,\sqrt{-1}\omega)$, where $(g,\omega)$ is the Weierstrass
data of $f$.
\begin{corollary}[A duality between swallowtails and cuspidal cross caps]
 \label{cor:desitter-duality}
 A \cmcone{} face $f:U\to S^3_1$ is $\A$-equivalent to a 
 swallowtail {\rm (}resp.\ a cuspidal cross cap{\rm )} at $p\in U$
 if and only if its conjugate \cmcone{} face $\tilde f$
 is $\A$-equivalent to a cuspidal cross cap 
 {\rm (}resp.\ a swallowtail{\rm )} at $p\in U$.
\end{corollary}



\begin{thebibliography}{KRSUY}
\bibitem[AGV]{AGV}
    V. I. Arnol'd, S. M. Gusein-Zade and A. N. Varchenko,
    {\itshape Singularities of differentiable maps, Vol. $1$},
    Monographs in Math. {\bf 82}, Birkh\"auser (1985).
\bibitem[Br]{Br}
        R. Bryant,
        {\itshape Surfaces of Mean Curvature One in Hyperbolic Space},
        Ast\'{e}risque {\bf 154-155} (1987), 321--347.
\bibitem[C]{C}
    J. P. Cleave,
    {\itshape The form of the tangent developable at points 
         of zero torsion on space curves},
    Math. Proc. Camb. Phil. Soc. {\bf 88}, (1980) 403--407.
\bibitem[FLS]{FLS}
    I. Fern\'andez, F. J. L\'opez and R. Souam, 
       {\itshape
         The space of complete embedded maximal surfaces
         with isolated singularities in the $3$-dimensional
         Lorentz-Minkowski space $\mathbb{L}^3$,
         }
        Math.\ Ann.\ {\bf 332} (2005), 605--643.
%
\bibitem[F]{F}
    S. Fujimori,
       {\itshape Spacelike CMC $1$ surfaces with
            elliptic ends in de Sitter $3$-Space},
        Hokkaido Math. J., {\bfseries 35} (2006), 289--320.
%
\bibitem[FRUYY]{FRUYY}
  S. Fujimori, W. Rossman, M. Umehara, K. Yamada and S.-D. Yang,
  {\itshape Spacelike mean curvature one surfaces in de Sitter
	    $3$-space}, preprint, arXiv:0706.0973.
%
\bibitem[GG]{GG}
    M. Golubitsky and V. Guillemin,
    {\itshape Stable Mappings and Their Singularities},
    Graduate Texts in Math. {\bf 14}, Springer-Verlag (1973).
%
\bibitem[IM]{IM}
    G. Ishikawa and Y. Machida,
        {\itshape Singularities of improper affine spheres and 
         surfaces of constant Gaussian curvature},
    International J. of Math. {\bf 17} (2006), 269--293.
%
\bibitem[ISTa]{ISTa}
  S. Izumiya, K. Saji and M. Takahashi,
  {\itshape Horospherical flat surfaces in hyperbolic $3$-space},
  preprint,
  http://eprints.math.sci.hokudai.ac.jp/archive/00001688/.
%
\bibitem[ISTe]{ISTe}
        S. Izumiya, K. Saji and N. Takeuchi,
        {\itshape Circular surfaces},
    Advances in Geometry. {\bf 7} (2007), 295--313.
\bibitem[K]{K}
  O. Kobayashi,
  {\itshape
    Maximal surfaces with conelike singularities},
    J. Math.\ Soc.\ Japan {\bfseries 36} (1984), 609--617.
\bibitem[KRSUY]{KRSUY}
    M. Kokubu, W. Rossman, K. Saji, M. Umehara and K. Yamada,
    {\itshape Singularities of flat fronts in hyperbolic $3$-space},
    Pacific J. Math. {\bfseries 221} (2005), 303--351.
\bibitem[SUY]{SUY} 
K. Saji, M. Umehara and K. Yamada,
        {\itshape The geometry of fronts}, 
         to appear in Ann. of Math., math.DG/0503236.
\bibitem[UY1]{UY1}
        M. Umehara and K. Yamada,
        {\itshape Complete surfaces of constant mean curvature $1$ 
         in the hyperbolic $3$-space},
         Ann. of Math. (2) {\bf 137} (1993), 611--638.
\bibitem[UY2]{UY}
        \bysame,
        {\itshape Maximal surfaces with singularities in 
         Minkowski space}, 
        Hokkaido Math. J., {\bfseries 35} (2006), 13--40.
\bibitem[W]{Whitney}
        H. Whitney, 
        {\itshape The singularities of a smooth $n$-manifold in
       $(2n-1)$-space},
      Ann.\ of Math.\ {\bfseries 45} (1944), 247--293.
\end{thebibliography}
\end{document}